\documentstyle[12pt,amssymb]{amsart}
\title
[Filling-Invariants at Infinity for Manifolds of Curvature${}\le0$]
{Filling-Invariants at Infinity for Manifolds of Nonpositive Curvature}
\author{Noel Brady}
\address{Noel Brady\\
Math Department, University of Utah\\
Salt Lake City, Utah 84112}
\email{brady@@math.utah.edu}
\author{Benson Farb}
\address{Benson Farb\\
Math Department, University of Chicago\\
5734 University Ave.\\
Chicago, Illinois 60637}
\email{farb@@math.uchicago.edu}
\thanks{Farb is supported in part 
by an NSF Postdoctoral Fellowship. His work at MSRI is supported by
NSF grant DMS-9022140.}

\newtheorem{theorem}{Theorem}[section]
\newtheorem{prop}[theorem]{Proposition}
\newtheorem{lemma}[theorem]{Lemma}

\def\proof{{\bf {\medskip}{\noindent}Proof: }}

\def\examples{{\bf {\bigskip}{\noindent}Examples: }}
\def\question{{\bf {\bigskip}{\noindent}Question: }}
\def\dfinition{{\bf {\bigskip}{\noindent}Definition: }}
\def\remark{{\bf {\bigskip}{\noindent}Remark: }}

\def\eproof{$\Box$ \bigskip}

\def\epar{\medskip}
\def\df{\em}
\def\title{\em}

\def\sdirect{\times \nobreak\kern-.19em\rule{.25pt}{1.1ex}}
\def\rtimes{\times \nobreak\kern-.19em\rule{.25pt}{1.1ex}}

\def\H{\mbox{\bf H}}
\def\R{\mbox{\bf R}}

\def\Sol{\mbox{\bf Sol}}
\def\div{\mbox{div}}

\begin{document}
\input{psfig}

\maketitle

\setcounter{section}{-1}

\section{Introduction}

Homological invariants ``at infinity'' and (coarse) 
isoperimetric inequalities are basic tools in the study of 
large-scale geometry (see e.g., \cite{Gr}).  
The purpose of this paper is to combine these two ideas to construct a
family $\div_k(X^n), \ 0\leq k\leq n-2$, of geometric invariants for 
Hadamard manifolds $X^n$
\footnote{Recall that a {\df Hadamard manifold} is a 
complete, simply-connected manifold with nonpositive sectional curvatures.}
.  The $\div_k(X^n)$ are meant to give a finer measure of 
the spread of geodesics in $X^n$; in fact the 
0-th invariant $\div_0(X^n)$ is the well-known 
``rate of divergence of geodesics'' in the Riemannian manifold $X^n$.

The definition of $\div_k(X^n)$ goes roughly as follows (see 
Section \ref{def} for the precise definitions): Find the minimum 
volume of a ball $B^{k+1}$ needed to fill a sphere $S^k$, 
where $S^k$ sits on the sphere $S(r)$ of radius $r$ in 
$X^n$, and the filling ball $B^{k+1}$ is required to lie outside the 
open ball $B(r)^{\circ}$ in $X^n$.  Then 
$\div_k(X^n)$ measures the growth of this volume as $r\rightarrow 
\infty$; hence $\div_k(X^n)$ is in some sense a $k$-dimensional 
isoperimetric function at infinity. 

We view the invariants $div_k(X^n)$ in the same way as we 
view the standard isoperimetric inequalities (for manifolds or for groups): 
as basic geometric quantities to be computed.

The $\div_k({X^n})$ are quasi-isometry invariants of $X^n$.  
The fundamental group $\pi_1(M^n)$ (endowed with the word metric) 
of a compact Riemannian manifold is quasi-isometric to the 
universal cover $\widetilde{M^n}$; hence the $div_k(\widetilde{M^n})$ 
give quasi-isometry invariants for fundamental groups of closed, 
nonpositively curved manifolds $M^n$.

The contents of this paper are as follows:  In Section \ref{def}, 
$div_k(X^n)$ is defined and shown to be a quasi-isometry invariant.  
The core of this paper (Sections \ref{suspend 
section},\ref{pulloff section},\ref{hyper section}) 
describes three geometric techniques for computing $div_k(X^n)$ for some 
basic examples.  
Section \ref{hyper section} also explores some surprising 
quasi-isometric embeddings hyperbolic spaces and solvable Lie groups 
into products of hyperbolic spaces.  
\medskip

We would like to express our gratitude to K. Fujiwara, C. Pugh, 
and A. Wilkinson for useful discussions, and to G. Kuperberg for doing 
the first two figures.

\section{Definitions and Quasi-isometry Invariance}
\label{def}

Let $X^n$ be a Hadamard manifold; that is,  a complete, simply connected 
manifold all of whose sectional curvatures are nonpositive. 
Let $S(r)$, $B(r)$ and $B(r)^\circ$ denote respectively the sphere, ball and 
open ball of radius $r$ about a fixed basepoint $x_0$ of $X^n$, and let   
$C(r)\, = \, X^n \setminus B(r)^\circ$. Note that $C(r)$ deformation 
retracts onto the sphere $S(r)$; hence any continuous map $f:S^{k} \to S(r)$
admits a continuous extension, or {\em filling} 
$\widehat{f}:B^{k+1} \to C(r)$, for any integer $0 \leq 
k \leq n-2$. 

We shall be considering lipschitz maps to the manifold $X^n$. 
By the Whitney Extension Theorem, we know that if 
$f$ as above is lipschitz,  then the extension 
$\widehat{f}$ of $f$ can be chosen to be lipschitz (with the same 
lipschitz constant). 
By Rademacher's Theorem, lipschitz maps are differentiable almost everywhere, 
enabling one to define the {\em $k$-volume} of $f:S^k \to X^n$ and the 
{\em $(k+1)$-volume} of $\widehat{f}:B^{k+1} \to X^n$, 
where $S^k$ and $B^{k+1}$ denote the 
unit sphere and ball in euclidean space $\R^n$. 
More precisely, if the derivative $D_xf$ exists 
at a point $x \in S^k$, it sends an orthonormal basis at $x$ to a $k$-tuple of 
vectors in $T_{f(x)}(X^n)$. We can comupte the $k$-volume of the parallelopiped 
spanned by this $k$-tuple using the metric on $X^n$. This defines a function 
$V(x)$ almost everywhere on $S^k$, and we can then define the {\em $k$-volume} 
of $f$, denoted $\mbox{\rm vol}_k(f)$, to be the integral of $V$ over $S^k$.
This integral exists because $V(x)$ is a bounded 
measurable function defined almost everywhere on $S^k$, as $\|D_xf\|$ 
is bounded 
by the lipschitz constant of $f$. 

We are now ready to define the invariant $\div_k(X^n)$ for a fixed 
integer $0 \leq k \leq n-2$. Although the concept of $\div_k$ is quite simple, 
the precise definition of $\div_k$ needs to be somewhat technical in order to 
make it manifestly a quasi-isometry invariant. 
This is accomplished using a variation of a 
trick introduced in \cite{Ge}. 

Let $A > 0$ and $0 < \rho \leq 1$ be given. 
For $r > 0$, we define a map $f:S^k \to S(r)$ to be {\em $A$-admissible} if: 
\begin{itemize}
\item $f$ is lipschitz, and 
\item $\mbox{\rm vol}_k(f) \leq Ar^k$
\end{itemize}
and say that the extension $\widehat{f}$ of $f$ is {\em $\rho$-admissible} if:
\begin{itemize}
\item $\widehat{f}$ is lipschitz, and
\item $\widehat{f}(B^{k+1}) \subset C(\rho r)$. 
\end{itemize}

In other words, the only admissible fillings are those which lie outside the 
open ball $B^\circ(\rho r)$ in $X^n$. 
Now define
 
\[\delta^k_{\rho, A}(r) \; = \; \sup_f\inf_{\widehat f}\mbox{\rm vol}_{k+1}(
\widehat{f}) \]

\noindent
where the supremun and infimum are taken over $A$-admissible maps $f$ and 
$\rho$-admissible fillings $\widehat f$ of $f$. 
We call the resulting two-parameter family of functions 

\[\div_k(X^n) \; = \; \{\delta^k_{\rho, A} \; : \; 
0 < \rho \leq 1, \; A>0 \}\]

\noindent
the {\em $k$-th divergence} of $X^n$ with respect to the 
point $x_0$. The parameters $\rho$ and $A$ are necessary in order to make 
$\mbox{\rm div}_k$ into a quasi-isometry invariant (see 
Theorem~\ref{invariance}). 

\remark We note that the function $\delta^k_{\rho, A}(r)$, as a $\sup$ of an 
$\inf$, may not be realized by an 
actual filling, though of course there 
are (admissible) fillings arbitrarily 
close to realizing this function.  We will 
ignore this distinction in what follows, as we are only interested 
in the {\em growth} of $\delta^k_{\rho, A}(r)$.
\medskip

In this paper we shall only be concerned with distinguishing between 
polynomial and exponential functions. Hence the following 
equivalence relation:  given functions $f,g:\R^+\rightarrow \R^+$, we write 
$f \preceq g$ if there exist constants $a, b, c>0$ and an 
integer $s\geq 0$ such that 
$f(x) \leq ag(bx) +cx^s$ for all sufficiently large $x$.  
Now write $f \sim g$ if both 
$f \preceq g$ and $g \preceq f$. 
This defines an equivalence relation on the class of functions from $\R^+ 
\to \R^+$, and it makes sense to call the equivalence classes 
polynomial, exponential, super-exponential, etc. 

Similarly, one defines an equivalence 
relation among $k$-th divergences as follows: say that $\mbox{\rm div}_k 
\preceq \mbox{\rm div}'_k$ if there exist $0 < \rho_0, \rho'_0 \leq 1$ and 
$A_0, A'_0>0$ such that for every pair $(\rho, A)$ with $\rho < \rho_0$ and 
$A > A_0$ there exist $\rho' < \rho'_0$ and $A' > A'_0$ with $\delta^k_{\rho, 
A} \preceq {\delta'}^k_{\rho', A'}$. 
Now define $\mbox{\rm div}_k \sim \mbox{\rm div}'_k$ if we have both 
$\mbox{\rm div}_k \preceq \mbox{\rm div}'_k$ and $\mbox{\rm div}'_k \preceq 
\mbox{\rm div}_k$. In particular, we say that $\div_k$ is polynomial or  
exponential, written $\div_k = poly$ or $\div_k = exp$, if there exists 
$0 < \rho_0 \leq 1$ and $0 < A_0$ such that $\delta_{\rho, A}(r) \sim r^d$ 
(for some integer $d>0$)  
or $\delta_{\rho, A}(r) \sim e^r$ for all $\rho<\rho_0$ and $A > A_0$. 
Thus one can speak of polynomial or exponential $k$-th divergences. 
\bigskip

We are now ready to prove that the invariants $div_k(X^n)$ are 
actually quasi-isometry invariants of $X^n$, sometimes called 
{\em geometric invariants}.  Recall that 
a quasi-isometry is basically a coarse bi-lipshitz map; these are the 
appropriate maps to study when one is interested in large-scale geometric 
properties of a space, or in geometric properties of the fundamental group 
of a compact Riemannian manifold (see, e.g., \cite{Gr}).  More precisely, 
we recall the following:

\dfinition Let $X$ and $Y$ be metric spaces.   
A {\em quasi-isometry} is a pair of maps 
$f:X\longrightarrow Y, g:Y\longrightarrow X$ 
such that, for some constants $K,\epsilon,C>0$ :

\[
\begin{array}{ll}
d_Y(f(x_1),f(x_2))\leq Kd_X(x_1,x_2)+\epsilon, &
d_X(g\circ f(x_1),x_1)\leq C\\
d_X(g(y_1),g(y_2))\leq 
Kd_Y(y_1,y_2)+\epsilon, &
d_Y(f\circ g(y_1),y_1)\leq C
\end{array}
\]

\noindent
for all $x_1,x_2\in X, y_1,y_2\in Y$.  Note that 
neither $f$ nor $g$ need be continuous.  If such maps exist, 
$X$ and $Y$ are said to be {\em quasi-isometric}; 
the map $f$ is called a $K$-quasi-isometry.  A {\em quasi-isometric 
embedding} is defined similarly.  A basic example to keep in mind is 
that the fundamental group $\pi_1(M)$ (endowed with the word metric) 
of a compact Riemannian manifold 
$M$ is quasi-isometric to the universal cover $\widetilde{M}$ of $M$.  
\medskip

\begin{theorem}[$div_k$ is a quasi-isometry invariant]
\label{invariance}
The $k$-th divergence, $\mbox{\rm div}_k$, is a quasi-isometry invariant of 
Hadamard manifolds.  In particular, $\mbox{\rm div}_k$ gives a 
quasi-isometry invariant for fundamental groups of closed, 
nonpositively curved manifolds.
\end{theorem}

Theorem \ref{invariance} allows us to simply speak of 
the $k$-th divergence $div_k(X^n)$ as ``polynomial'' or ``exponential'', 
denoted by $\div_k(X^n)\sim exp$ and $\div_k(X^n)\sim poly$, 
respectively, without 
having to speak of the actual 2-parameter family of functions by which 
$div_k(X^n)$ is defined.  Whether $div_k(X^n)$ is polynomial or 
exponential is then a quasi-isometry invariant notion.

\bigskip
\noindent
{\bf Proof of Theorem \ref{invariance}:  }
Let $F:(X,x_0) \to (Y,y_0)$ and $G:(Y,y_0) \to (X,x_0)$ be 
$L$-lipschitz maps between two based Hadamard manifolds 
which are determined by a quasi-isometry  between $X$ and 
$Y$ (see Appendix~A for existence of lipschitz quasi-isometries). We 
shall use these maps to compare the $k$-th divergences $\div_k$ for 
$X$ and $\div'_k$ for $Y$. 

\begin{figure}[th]
\label{divk figure}
\centerline{\psfig{file=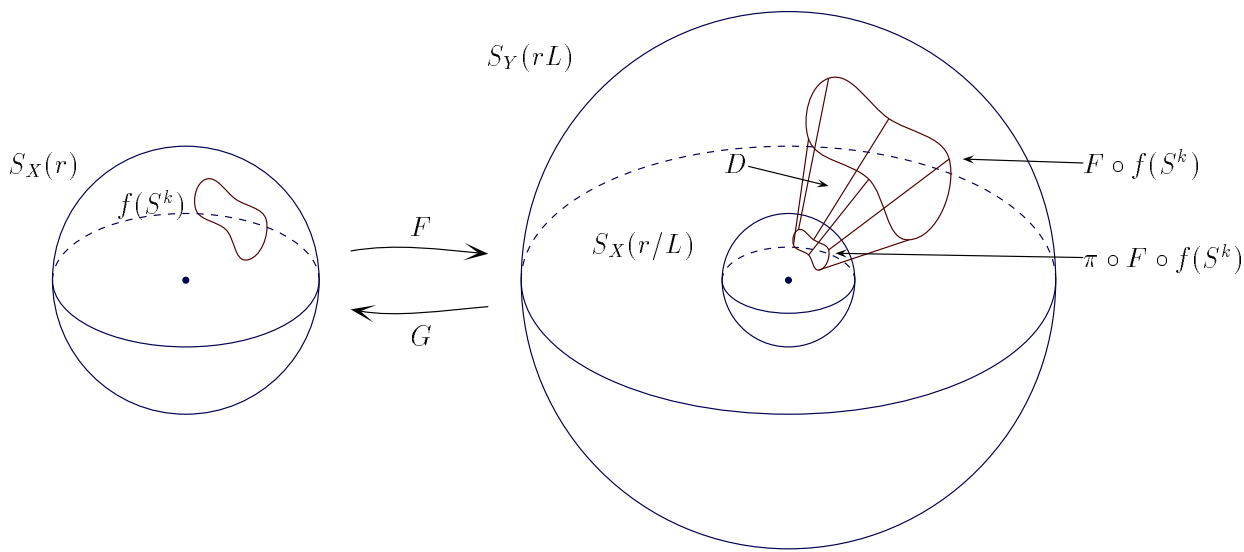,height=231pt,width=470pt}}
\caption{This figure illustrates part of the proof of Theorem 
\protect \ref{invariance}.}
\end{figure}

Given an $A$-admissible map $f:S^k \to 
S_X(r)$, we can compose with $F$ to get a lipschitz map $F\circ f: 
S^k \to Y$ with $(F\circ f)(S^k) \subset B_Y(Lr) \cap C_Y(r/L)$. Note that 
$\mbox{\rm vol}_k(F\circ f) \leq A(Lr)^k$. Radial projection onto $S_Y(r/L)$ 
defines a volume non-increasing lipschitz map (lipschitz constant 1) 
$\pi: C_Y(r/L) \to S_Y(r/L)$, and so the composition $(\pi \circ F\circ f): 
S^k \to S_Y(r/L)$ is $AL^{2k}$-admissible (see Figure 
\ref{divk figure}). There is an admissible  
filling of this map with $(k+1)$-volume 
bounded by $\delta^{\prime k}_{\rho , A/L^{2k}}(r/L)$.

We obtain a lipschitz extension of 
$F\circ f$ as follows: on the radius one-half ball in $B^{k+1}$ take the map 
$(\widehat{\pi\circ F\circ f})\circ d_2$ where $d_2$ denotes the dilation 
taking the radius one-half ball $\frac{1}{2}B^{k+1}$ 
onto $B^{k+1}$, and on the remaining 
annular region just interpolate between the maps $F\circ f$ and 
$(\widehat{\pi \circ F\circ f})\circ d_2$ (sending radial geodesics in 
$B^{k+1}$ to geodesic segments between image points in $Y$). Note that the 
$(k+1)$-volume of this map is bounded by $\delta'_{\rho, AL^{2k}}(r/L)$ plus a 
polynomial of degree $k+1$ in $r$. This polynomial bounds the $(k+1)$-volume 
of the annular region, and is the reason we use the equivalence relation 
among functions defined above.

Now postcomposition with $G$ yields a lipschitz map from $B^{k+1}$ to 
$X$ which lies outside the $\rho r/l^2$-ball about $x_0$, and the restriction 
of this map to $S^k$ is a constant distance (pointwise) away from the original 
map $f:S^k \to X$, as $G\circ F$ is a constant distance away from the 
identity map $1_X$.  It is easy to see that one can interpolate 
between these maps to obtain a lipschitz map $\widehat{f}:B^{k+1} \to X$ which 
is a $\rho/L^2$-admissible filling of $f$. Note that 

\[\mbox{vol}_{k+1}(\widehat f) \, \leq \, L^{k+1}\delta'_{\rho, AL^{2k}}(r/L) 
\, + \, p(r) \]

\noindent
where $p(r)$ is a polynomial of degree $k+1$ in $r$, and so 

\[\delta_{\rho/L^2, A}(r) \preceq \delta'_{\rho, AL^{2k}}(r) \, .\]

\noindent
Thus, $\div_k \preceq \div'_k$. Similarly, $\div'_k \preceq \div_k$ and 
so $\div_k$ is a true quasi-isometry invariant.

\eproof

The proof of Theorem~\ref{invariance} shows that $\div_k$ can be made more 
precise than $poly$ or $exp$; in fact, $\div_k$ is a well-defined 
quasi-isometry invariant up to an additive factor of $x^{k+1}$.  

\section{Suspending Hard-to-Fill Spheres}
\label{suspend section}

In this section we show how to suspend hard-to-fill spheres in $X$ to 
hard-to-fill spheres in $X\times \R$.  This provides a lower bound for 
the $(k+1)$st-divergence of $X\times \R$ in terms of the $k$th-divergence 
for $X$.
\medskip

\begin{theorem}[suspending hard-to-fill spheres]
\label{theorem1}
Let $(X^n,x_0)$ be a based Hadamard manifold.  Then 

\[\div_{k+1}(X^n\times \R)\geq \div_{k}(X^n) \mbox{\ \ \  for 
any\ \ }0\leq k\leq n-3.\]
\end{theorem}

\medskip
Theorem \ref{theorem1} shows, for example, that 

\[\div_2(\H^2\times \R^2)\geq \div_1(\H^2\times \R)\geq \div_0(\H^2)
\sim exp.\]

\begin{figure}[h]
\label{suspend fig}
\begin{center}~
\psfig{file=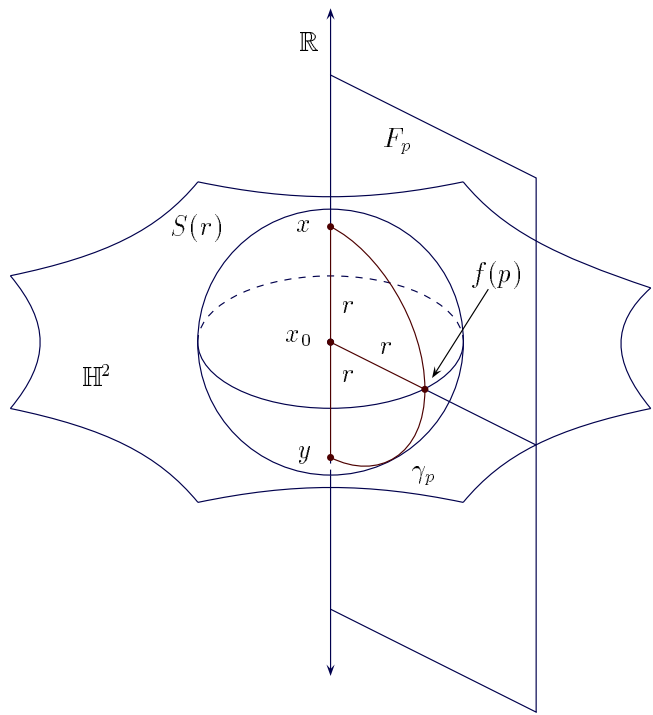,height=280pt,width=340pt}
\caption{The arc $\gamma_p$ is used to define $\Sigma(f)$ at the point 
$\{p\}\times [0,1]$ of the suspension $\Sigma(S^k)=S^{k+1}$.}
\end{center}
\end{figure}

{\bf {\medskip}{\noindent}Proof of Theorem \ref{theorem1}: } 
Denote $X^n$ simply by $X$, and let $S(r)$ be the sphere of radius $r$ in 
$X\times \R$.  Since $X=X\times\{0\}$ is totally geodesic in 
$X\times \R$, the intersection $S^{\prime}(r)=S(r)\cap X$ is the sphere 
of radius $r$ in $X$.  Choose an admissible map 
$f:S^k\rightarrow S^{\prime}(r)$ which 
realizes $div_k(X)$.  We now define a map
\[\Sigma(f):\Sigma(S^k)=S^{k+1}\longrightarrow S(r)\]

\noindent
where $\Sigma(S^k)$ denotes the suspension 
$S^k\times [0,1]/((S^k\times\{0\})\cup (S^k\times\{1\}))$ of $S^k$.  
Geometrically, define the map $\Sigma(f)$ as follows: 
for $p\in S^k$, let $F_p$ be the 2-flat 
which is the product of the infinite geodesic in $X$ 
passing through $x_0$ and $f(p)$ and the infinite geodesic 
$\{x_0\}\times \R$.  We think of each $F_p$ based at the ``origin'' 
$(x_0,0)\in X\times \R$, and note that the points 
$x=(x_0,r)$ and $y=(x_0,-r)$ of $S(r)$ are contained in each 
$F_p$.  Let $\gamma_p$ denote the arc of the circle of radius 
$r$ in $F_p$ from $x$ to $y$; this arc has length $\pi r$ (see Figure 
\ref{suspend fig}).  

Define the map $\Sigma(f)$ via:

\[\Sigma(f)(\{p\}\times [0,1])=\gamma_p.\]

So, for example, $\Sigma(f)$ stretches each $\{p\}\times [0,1]$ by 
a factor of $\pi r$.  It is not difficult to check that 
$\Sigma(f):S^{k+1}\rightarrow S(r)$ is an admissible map.  Note also 
that $\Sigma(f)(S^{k+1})\cap (X\times\{0\})=f(S^k)$.

Suppose that $\widehat{\Sigma(f)}:B^{k+1}\rightarrow C(r)$ is an admissible 
filling of $\Sigma(f)$.  Then $\widehat{\Sigma(f)}\cap(X\times\{0\})$ 
gives an admissible filling of $f(S^k)$ in $X$, and this filling 
has volume at least $\delta^k(r)$, where $\delta^k(r)$ is the appropriate 
divergence function as in the definition of $div_k(X^n)$.

Let a constant $\tau <<1$ be fixed.  The reasoning above shows that for 
any $0\leq \epsilon \leq \tau$, the map given by 
$\widehat{\Sigma(f)}\cap(X\times\{\epsilon\})$ 
gives an admissible filling of $f(S^k)$ in $X\times \{\epsilon\}$, 
and this filling 
has volume at least $\delta^k(r-\epsilon)$.  Since the leaves 
$\{X\times \{\epsilon\}\}_{0\leq \epsilon \leq \tau}$ are parallel, we 
have

\[
\begin{array}{ll}
vol(\widehat{\Sigma(f)})&\geq \int_0^{\tau}\delta^k (r-\epsilon)d\epsilon\\
&\\
&\geq \int_0^{\tau}\delta^k (r-\tau)d\epsilon\\
&\\
&\geq \tau \cdot \delta^k(r-\tau)\\
&\\
&\sim \delta^k(r)
\end{array}
\]

\noindent
as desired.

\eproof

\section{Pulling-Off Spheres Along Flats}
\label{pulloff section}

In this section we give a polynomial upper bound for $div_1(X\times \R^2)$ 
for any Hadamard manifold $X$ (Theorem \ref{theorem2}).  
In order to do this we need the 
following, easily believable technical lemma.  The main idea in the proof 
of Theorem \ref{theorem2} may be digested indpendently of the proof 
of this lemma.

\begin{lemma}[homotoping off a neighborhood of the origin]
\label{perturb}
Let $\beta:S^1\rightarrow \R^2$ be a lipschitz map such that the length 
of $\beta(S^1)$ is at most $Ar$.  Then it is possible to homotope 
$\beta$ to $\beta':S^1\rightarrow \R^2$ of length at most $\pi Ar$,  
so that $\beta'(S^1)$ lies 
outside the open unit ball in $\R^2$. The paths $\beta$ and $\beta'$ are 
homotopic by a lipschitz homotopy 
of area at most $2\pi A r$. 
\end{lemma}

\proof
Let $B^\circ$ and $\partial B$ denote respectively  the open unit ball 
and unit circle in $\R^2$. The map $\beta$ is 
lipschitz and therefore continuous, and so $\beta^{-1}(B^\circ)$ is an open 
subset of $S^1$. There are three cases to consider. 

In the first case $\beta^{-1}(B^\circ) = \emptyset$, and we take $\beta' = 
\beta$ and the result is trivial. 

In the second case $\beta^{-1}(B^\circ) = S^1$. Here we take $\beta' = 
\tau_{\vec{v}} \circ \beta$ where $\tau_{\vec{v}} : \R^2 \to \R^2$ is just 
translation by a vector $\vec{v}$  of length 2. The homotopy is 
given by maps $\tau_{t \vec{v}}\circ \beta$ 
where $t \in [0,1]$. Since the $\tau_{t\vec{v}}$ are isometries of $\R^2$ 
in the Euclidean metric, $\beta'$ is lipschitz of length $Ar$, and the 
homotopy is clearly lipschitz, of area 
at most $2Ar$.

Finally, $\beta^{-1}(S^1)$ may be a non-empty proper open subset of $S^1$, 
and so is the disjoint union of a collection of open intervals in $S^1$. Given 
such an open interval $(a,b) \subset S^1$, we may define $\beta'|_{[a,b]}$ 
to agree with $\beta$ on the endpoints $a$ and $b$, and to map $[a,b]$ 
uniformly over the smaller of the arcs of $\partial B$ determined by 
$\beta(a)$ and $\beta(b)$ (either arc if $\beta(a)$ and $\beta(b)$ are 
antipodal).  We take the straight line homotopy $t\beta'(x) + (1-t)\beta(x)$ 
in $\R^2$, between $\beta$ and $\beta'$. It is clear from the construction
 that $\beta'$ and the homotopy are lipschitz, and that the length of $\beta'$ 
is bounded by $\pi Ar$, and that the area of the homotopy is 
at most $2\pi Ar$. 
\eproof

\begin{theorem}[a polynomial filling]
\label{theorem2}
Let $X$ be a Hadamard manifold.  Then 
$$
\div_1(X\times \R^2)
$$ is 
polynomial of degree three.
\end{theorem}

\proof
Let an $A$-admissible map $\gamma:S^1\rightarrow S(r)$ be given, where 
$S(r)$ is the sphere of radius $r$ around a chosen basepoint 
$x_0\in X$.  Let 
$\pi:X\times \R^2\rightarrow \R^2$ be the natural projection.  

From Lemma \ref{perturb} applied to $\beta=\pi\circ \gamma$, it is clear 
that we may homotope $\gamma$ slightly so that 
$\pi\circ\gamma(S^1)$ lies outside the open 
unit ball in $\R^2$.  An admissible 
filling of this perturbed $\gamma$ then gives an admissible 
filling of $\gamma$.  Since the 
area of these two maps differs by at most some constant (not depending 
on $r$) times $r$, we may assume without loss of generality 
that $\pi\circ\gamma(S^1)$ lies outside the open unit 
ball in $\R^2$.

\begin{figure}
\label{pulloff figure}
\begin{center}~
\psfig{file=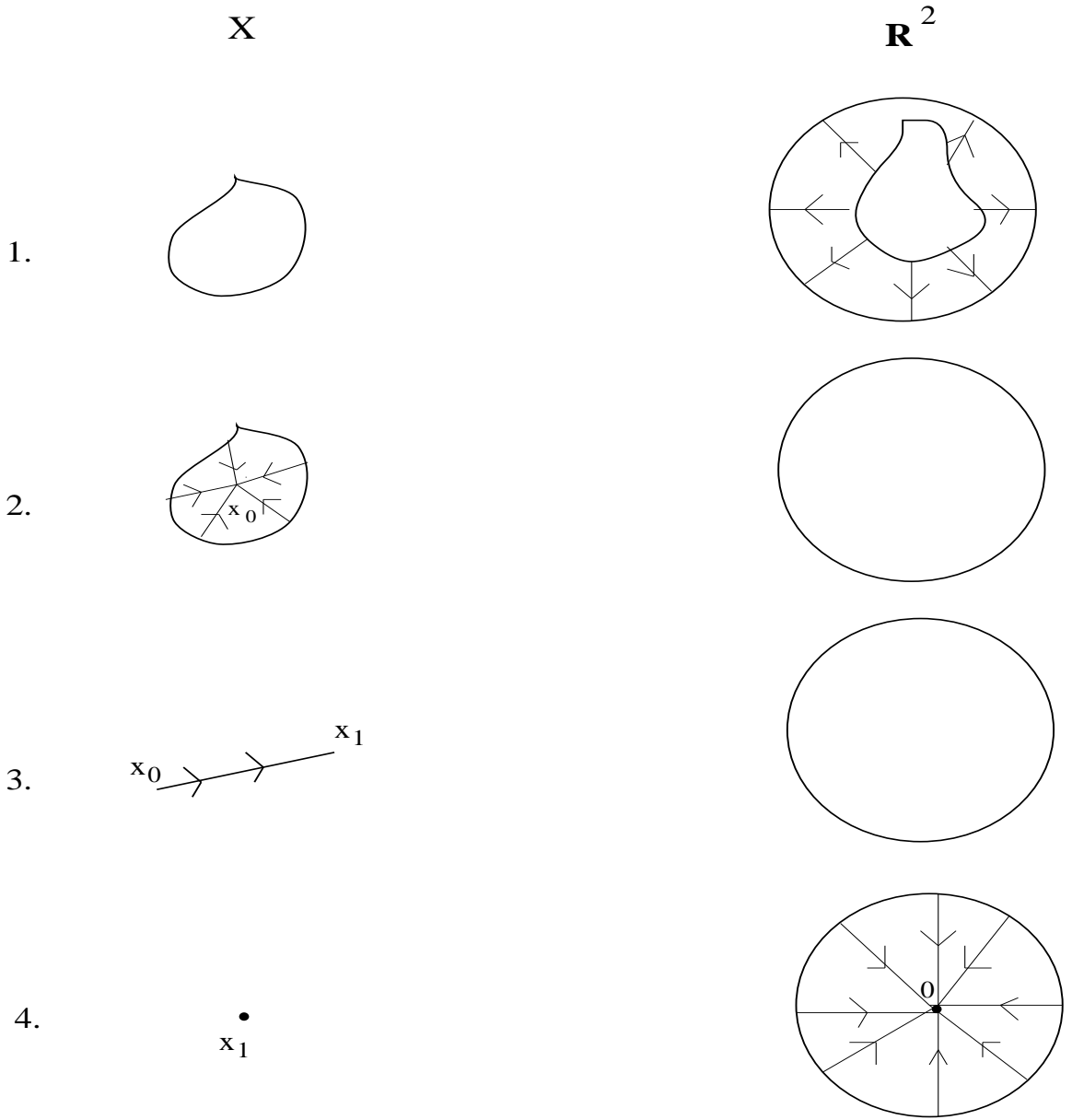,height=6.4in,width=6in}
\caption{This figure gives a sequence of 4 homotopies which homotope 
$\gamma$ (outside $B(r)$) to a point.  The homotopy is illustrated 
by its projections onto the $X$ and $\R^2$ factors.}
\end{center}

\end{figure}

We now give an admissible filling $\hat{\gamma}:B^2\rightarrow C(r)$ 
of $\gamma$.  The filling is given by the tracks of $\gamma$ under a 
sequence of 4 homotopies which homotope $\gamma$ (outside $B(r)$) 
to a point; this filling is illustrated in Figure 
\ref{pulloff figure}.

Here is the sequence of homotopies :

\begin{enumerate}

\item Radial projection in the $\R^2$ factor to the $3r$-sphere in 
$\R^2$: 
$(x,\vec{v})\mapsto (x,t\cdot\frac{\vec{v}}{||\vec{v}||}), 1\leq t\leq 3r$.  
Since $\pi\circ \gamma$ lies outside 
the open unit ball in $\R^2$, this radial projection is 
well-defined, and increases the length of $\gamma$ by at most a 
factor of $3r$ (hence the length of the image of 
$\gamma$ under this radial projection is at most $3r\cdot |\gamma|\leq 
3Ar^2$).

\item Coning-off in the $X$ factor: $(x,\vec{v})\mapsto 
(\sigma_x(t),\vec{v}), 0\leq t\leq 1$, where $\sigma_x$ is the (unique) 
geodesic in $X$ from $x$ to $x_0$.

\item Pulling to the $5r$-sphere in the $X$ factor: 
$(x_0,\vec{v})\mapsto (\tau(t),\vec{v}), 0\leq t\leq 1$, 
where $x_1$ is any (fixed) 
point lying on the sphere of radius $5r$ in $X$, and $\tau$ is the unique 
geodesic from $x_0$ to $x_1$.   

\item Coning-off in the $\R^2$ factor: 
$(x_1,\vec{v})\mapsto (x_1,(1-t)\cdot\vec{v}),0\leq t\leq 1$.

\end{enumerate}

It is easy to check that the images of these homotopies 
lies outside $B(r)$, since the metric on $X\times \R^2$ is 
just the product metric.  Piecing together these homotopies gives 
a map $\gamma':S^1\times [0,1]\rightarrow C(r)$ with the image 
$\gamma'(S^1\times \{1\})$ being the point $(x_1,\vec{0})$; hence 
this induces a map on the cone on $S^1$, that is a map 
$\hat{\gamma}:B^2\rightarrow C(r)$.  The map $\hat{\gamma}$ is 
easily seen to be an admissible filling of area at most

\[(3Ar^2)\cdot (3r)+(Ar)\cdot r+(3Ar^2)\cdot (5r)+(3Ar^2)\cdot (3r)
\leq 35Ar^3.\]

Hence $div_k(X\times \R^2)$ is (equivalent to) a polynomial of degree 
$3$.
\eproof

\remark It is probably true, more generally, that $div_k(X\times \R^m)$ is 
polynomial for $k<m$.  The proof of Theorem \ref{theorem2} works 
verbatim in this case, except for Lemma \ref{perturb} (which we 
believe to be true, although we have not been able to find a proof).

\section{Transverse Flats and Hyperbolic Spaces in Products}
\label{hyper section}

It is a surprising (but not difficult to see) fact that there 
is, for example, a quasi-isometrically embedded copy of 
$\H^3$ inside $\H^2\times \H^2$ 
\footnote{We somehow remember that this fact was stated by Gromov 
somewhere in \cite{Gr}, but we've been unable to locate the exact 
reference.}.  
In the first part of this 
section, we show that there are quasi-isometric embeddings of 
hyperbolic spaces and solvable Lie groups in 
products of hyperbolic spaces, although these embeddings are not 
quasi-convex.  We then apply the first embeddings 
to find a lower bound for certain 
$\div_k$ of products of hyperbolic spaces.  

The idea is to exploit the 
fact that, in a product of hyperbolic spaces $X$, there is a flat and 
a nicely embedded hyperbolic space whose dimensions add to $dim(X)+1$.  
The flat is used to find a polynomial-volume sphere; the hyperbolic 
space is used to show that any admissible filling of this sphere has 
exponential volume.

Recall that a subset $Y$ of a metric space $X$ is called {\em quasiconvex} 
in $X$ if there is a constant $D\geq 0$ so that any geodesic in $X$ 
between points $y_1,y_2\in Y$ lies in the $D$-neighborhood of $Y$.

\begin{prop}[q.i. embeddings]
\label{embedprop}
There are quasi-isometric embeddings of 
$$
\H^{(m_1\,+\,\cdots\, +\,m_n)\,-\,k+1}
$$ and of an $((m_1+\cdots +m_n)-k+1)$-dimensional 
solvable Lie group in 
$X=\H^{m_1}\times \cdots \times \H^{m_k}$, where each $m_i>1$.  
These embeddings are not quasiconvex.
\end{prop}

\examples Proposition \ref{embedprop} shows that there are quasi-isometrically 
embedded (but not quasiconvex) copies of $\H^3$ and of the three-dimensional 
geometry $\Sol$ in $\H^2\times \H^2$, and that there are quasi-isometrically  
embedded copies of $\H^4$ in $\H^2\times\H^2\times\H^2$ and of 
$\H^5$ in $\H^3\times \H^3$.
\medskip

\proof
Let $\gamma(t) = (\gamma_1(t), \ldots , \gamma_k(t))$ be an 
infinite diagonal geodesic in $X$, parameterized by arc length; so that 
$\gamma_i$ traces out (at a speed of $1/{\sqrt{k}}$ times unit speed) 
an infinite geodesic 
$g_i$ in $H^{m_i}$. For each $1\leq i \leq k$, let $S_i^+(t)$ (respectively 
$S_i^+(t)$) denote the 
horosphere in $\H^{m_i}$ which is centered at $\gamma_i(+\infty)$ 
(respectively $\gamma_i(-\infty)$) and which contains the point $\gamma(t)$. 
\epar

Now we define the quasi-isometrically embedded copy of hyperbolic space
\footnote{When using the coordinates $\R^{n-1}\sdirect\,\R$ for hyperbolic 
space $\H^n$ (where $t \in \R$ acts on $\R^{n-1}$ by  $(x_1, \ldots ,x_n) 
\mapsto e^{-t}(x_1, \ldots ,x_n)$) we shall refer to the $\R^{n-1}$ coordinates as {\em horospherical} and the $\R$ coordinate as {\em vertical}.}
$$Y \; = \; \H^{(m_1 + \cdots +m_k) -k+1} \; = \; \R^{(m_1 + \cdots + m_k) -k} 
\sdirect\,\R$$ 
in $\H^{m_1} \times \cdots \times \H^{m_k}$ to be the set of points 
$$\{ (x_1, \ldots , x_k) \; : \; x_i \in S_i^+(t) \; , \; t\in \R \}$$
Note that this set contains the geodesic $\gamma$ which is the 
vertical $\R$  in \\
$\R^{(m_1 + \cdots + m_k) -k} \sdirect\,\R$. 
Each horosphere $S_i$ contributes a copy of $\R^{m_i -1}$  which is scaled by 
a factor of $e^{-t/{\sqrt{k}}}$. Hence the metric induced on \\
$\R^{(m_1 + \cdots + m_k) -k} \sdirect\,\R$ is given by 
$$dt^2 \, \, + e^{-t/{\sqrt{k}}}ds^2$$
where $ds^2$ denotes the Euclidean metric on $\R^{(m_1 + \cdots + m_k)-k}$. 
\epar

\begin{figure}
\begin{center}~
\psfig{file=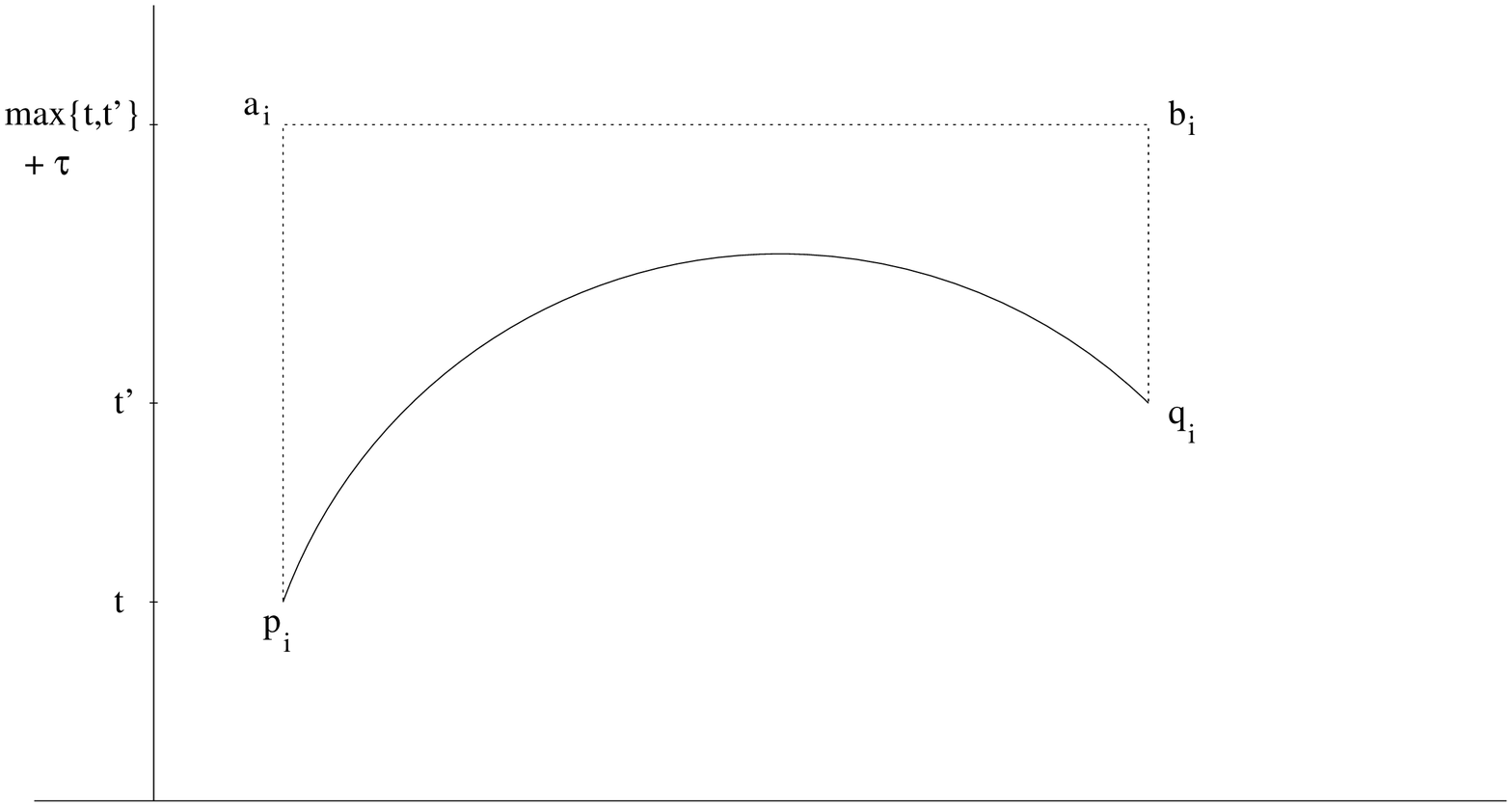,height=2.9in,width=3.4in}
\end{center}
\caption{The hyperbolic space is quasi-isometrically embedded.}
\label{qiembedding}
\end{figure}

First we show that $Y \subset X$ is a quasi-isometric embedding. For any two 
points $p,q \in Y$ we have $d_X(p,q) \leq d_Y(p,q)$. For the rest of 
the proof we refer to Figure \ref{qiembedding}. Let $p_i$ and 
$q_i$ denote the projections of $p$ and $q$ onto the factors $\H^{m_i}$, and 
let $\gamma_i$ denote the geodesic in $\H^{m_i}$ between $p_i$ and $q_i$. 
Since $p,q \in Y$, the $p_i$ all have the same vertical coordinate $t$ and 
the $q_i$ all have the same vertical coordinate $t'$. 
\epar

For each $1 \leq i \leq 
k$ let $$\pi_i:\H^{m_i} \to g_i : (x_1, \ldots , x_{m_i},t) \mapsto t$$
denote the horospherical projection onto the vertical geodesic $g_i$, and let 
$$\tau \, = \, \sup_i\mbox{length}(\pi_i(\gamma_i)) \, .$$ 
Denote by $a_i$ and $b_i$ respectively, the points in $\H^{m_i}$ which lie 
vertically above $p_i$ and $q_i$ on the horospherical level $\min\{t, t'\} \, 
+ \, \tau$.
\epar

Now consider the path in $Y$ from $p$ to $q$ which consists of the diagonal 
geodesic segments from $p = (p_i, \ldots ,p_k) $ to $(a_i, \ldots ,a_k)$ 
and from $(b_1, \ldots ,b_k)$ to $q = (q_1, \ldots ,q_k)$, and a geodesic 
(in the intinsic metric on the product of horospheres $S_1^+ \times \cdots 
\times S_k^+$) from $(a_1, \ldots , a_k)$ to $(b_1, \ldots ,b_k)$. 
Each of the diagonal geodesics have length bounded by $\sqrt{k}\tau$ and the 
other geodesic segment has length with bound $2\sqrt{k}$. 
\epar

Thus the length of this path is less than 
$2\sqrt{k} \tau + 2\sqrt{k}$ and so we have 
\begin{eqnarray*}
d_Y(p,q) &\leq & 2\sqrt{k}\tau + 2\sqrt{k} \\
 & \leq & 2\sqrt{k}d_{\H^{m_j}}(p_j,q_j) + 2\sqrt{k} \\
 & \leq & 2\sqrt{k}d_X(p,q) + 2\sqrt{k}
\end{eqnarray*}
where $\H^{m_j}$ is the factor realizing the maximum of the lengths of 
$\pi_i (\gamma_i)$. 
\epar


Now we show that $Y$ is not quasi convex in $X$. Consider the two points 
$(p_1, p_2, \ldots ,p_k)$ and $(q_1,p_2,\ldots,p_k)$ in $Y$ where $p_1$ 
and $q_1$ lie on the same horosphere in $\H^{m_1}$.  
The $X$ geodesic between these points is just the 
geodesic in $\H^{m_1}$ from $p_1$ to $q_1$ (with the other coordinates 
just constant at the point $(p_2, \ldots , p_k)$). Clearly this does not 
lie in $Y$; in fact, the distance from its midpoint $(r_1,p_2,\ldots,p_k)$ to 
$Y$ is given by 
$$\min_{l_1 + l_2 = l}\sqrt{(k-1)l_1^2 + l_2^2}$$
where $l$ is the vertical height of the geodesic between $p_1$ and $q_1$ 
in $\H^{m_1}$. This can be made arbitrarily large by choosing $p_1$ and 
$q_1$ far apart. 
\epar

Finally, note that there are quasi isometrically embedded copies of solvable 
Lie groups in $X$; these are just the horospheres in $X$. 
More explicitly, for example, define $Z$ in $X$ to be the set of points 
$$Z \; = \; \{(x_1, \ldots ,x_k)\; : \; x_1 \in S^+_1(t)\,,\;  x_i 
\in S^-_i(t)\; 
\; (i = 2, \ldots ,k)\, ,\;  t \in \R\}$$
We leave it to the reader to verify that the induced metric on 
$Z \; = \; \R^{(m_1\,  +\, \cdots \,+\,m_k) \, - \, k} \sdirect \, \R$ is 
given by $$dt^2 + \, e^{-t/\sqrt{k}}ds_1^2 + e^{t/\sqrt{k}}ds_2^2$$
where $ds_1^2$ denotes the Euclidean metric on $\R^{m_1 -1}$ and 
$ds_2^2$ denotes the Euclidean metric on $\R^{(m_2 + \cdots +m_k)-k+1}$, and 
that $Z$ is quasi isometrically embedded in $X$.
\eproof

The geodesic $\gamma$ of Proposition~\ref{embedprop} is contained in the 
$k$-flat $F=\gamma_1 \times \cdots \times \gamma_k$. This $k$-flat is foliated 
by parallel geodesics of the form 
$$\gamma_{(s_1, \ldots , s_k)}(t) \; = \; (\gamma_1(t+s_1), \ldots , 
\gamma_k(t+s_k))$$
where $s_1 + \cdots + s_k = 0$. The corresponding family of 
hyperbolic spaces
$$Y_{(s_1, \ldots , s_k)} \; = \; \{(x_1, \ldots , x_k) \; : \; x_i \in 
S^+_i(t+s_i) \, , \; \; t \in \R\}$$
gives a codimension $(k-1)$ 
foliation of $X$, where each $Y_{(s_i, \ldots , s_k)}$ intersects 
the $k$-flat $\gamma_1 \times \cdots \times \gamma_k$ in the geodesic 
$\gamma_{(s_1, \ldots, s_k)}$. 

\begin{lemma}
\label{stackofh}
For any points $p \in Y_{(s_1, \ldots, s_k)}$ and $q \in 
Y_{(s'_1, \ldots , s'_k)}$ we have 
$$d_X(p,q) \; \geq \; \sqrt{(s_1 -s'_1)^2 \, + \, \cdots \, + \, (s_k - s'_k)^2
}$$ 
\end{lemma}

\proof
Perpendicular projection of $X$ onto the $k$-flat $\gamma_1 \times \, 
\cdots \, \times \gamma_k$ is a distance nonincreasing map which takes 
$p$ and $q$ to points on the geodesic lines $\gamma_{(s_1, \ldots, s_k)}$ 
and $\gamma_{(s'_1, \ldots, s'_k)}$ respectively. The $X$-distance between the 
image points is the same as the distance in the $k$-flat, which is 
bounded below by $\sqrt{(s_1-s'_1)^2 \, + \, \cdots \, + \, (s_k - s'_k)^2}$. 
\eproof

\begin{theorem}[some hard-to-fill spheres]
\label{theorem3}
Let $X=\H^{m_1}\times \cdots \times \H^{m_k}$, each $m_i>1$, be a product of 
$k$ hyperbolic spaces.  Then $\div_{k-1}(X)=exp$.
\end{theorem}

\proof
Consider the family of hyperbolic spaces $\{Y_s\}$ as above, where 
$s\in B(1)$, the ball of radius 1 in $F$.  
Recall that $F=\gamma_1\times \cdots \times \gamma_k$ is 
a totally geodesic, isometrically embedded copy of $\R^k$ in $X$.  Also 
recall that the volume of a sphere of radius 
$r$ in $\R^k$ is $ar^k$ for some constant $a$ depending only on $k$.

Now let $A>a$ and $r>0$ be given, and 
let $T=F\cap S(r)$.  Then $T$ is a $(k-1)$-dimensional 
sphere lying on the sphere $S(r)$ of radius $r$ about the origin in $X$.  
Since $T\subset F$ and $F$ is a flat in $X$, we have 
$vol_{k-1}(T)\leq ar^k\leq Ar^k$.  We claim that 
any filling of $T$ outside of $S(r)$ has $(k+1)$-volume on the order 
of $e^r$.  

To prove this claim, 
suppose that $\widehat{T}$ is an admissible filling of $T$.  First 
note that the intersection $T\cap Y_s$ in $X$ has dimension 0; in fact, 
$T\cap Y_s$ consists precisely of the two points 
$T\cap \gamma_s=\{x_s,y_s\}$.  
Since $x_s$ and $y_s$ are each within a distance of 1 
from antipodal points of $S(R)$ lying on the 
geodesic $\gamma_0$, it follows that $x_s$ and $y_s$ are each within 
a distance of 1 from antipodal points 
on the sphere $S^{\prime}(r)$ of radius $r$ in $Y$.  Now 
$\widehat{T}\cap Y_s$ is a one-dimensional arc in $Y_s$ which connects 
$x_s$ to $y_s$ outside of $S^{\prime}(r)$.  Since $Y_s$ is a hyperbolic 
space, $\div_0(Y_s)$=exp, so that the arc $\widehat{T}\cap Y_s$ has length 
at least $Ce^r$ for some constant $C$ which is independent of $r$.  
We note that since $s\in B(1)$, the 
constant $C$ may be chosen to work for all $Y_s$.

Now

\[
\begin{array}{lll}
vol_{k+1}(\widehat{T})&=\int_{B(1)}vol_{k+1}(\widehat{T}\cap Y_s)d\mu(s)& \\
&&\\
&\geq \int_{B(1)}Ce^rd\mu(s)&\mbox{since $|\widehat{T}\cap Y_s|\geq Ce^r$}\\
&&\\
&\geq Ce^r\mu(B(1))&\mbox{by Lemma \ref{stackofh}}
\end{array}
\]

and we are done.  

\eproof

\remark The reason it is necessary to use 
Lemma \ref{stackofh} is that it is possible 
to have, for example, a disc foliated by an interval's worth of 
lines of length $e^r$, but with the area of the disc being constant.  
For example consider a long, thing quadrilateral in the hyperbolic plane: it 
is foliated by lines of length $e^r$ for $r$ large, but it's area is 
bounded by a universal constant; the reason is that the leaves of the 
foliation are bent so that they come very close together, on the order of 
$e^{-r}$, in fact.

\section{Questions}

As stated above, we view the invariants $div_k(X)$ as basic 
geometric quantities to be computed.  
We believe that the invariants $\div_k(X^n)$ are computable for many 
more examples than are covered in this paper.  

\question Compute $\div_k(X)$ for symmetric spaces $X$ of noncompact type.  
The simplest case not covered in this paper is $\div_2(\H^2\times \H^2)$.  
We believe that $div_1(X)\sim e^r$ for the symmetric space 
$X=SL_n(\R)/SO_n(\R)$.

\question Can the invariants $div_k(X)$ be used to detect the rank of a 
(globally) symmetric space $X$ of noncompact type?

\section*{Appendix A}

In this appendix we give a technical proposition which was needed 
to show that the $div_k(X^n)$ are quasi-isometry invariants.  
\bigskip

{\noindent \bf Proposition A.1: } 
{\em Suppose $f:X\rightarrow Y$ is a quasi-isometry 
between Hadamard manifolds, and suppose that $X$ admits a cocompact 
lattice.  Then $f$ is a bounded distance from a (continuous) lipshitz map 
$f':X\rightarrow Y$; that is, $sup_{x\in X}d_Y(f(x),f'(x))\leq C$ for 
some constant $C>0$.} 
\medskip

\proof
Since $X$ admits some compact quotient $M=X/{\Gamma}$, it is possible to lift 
a triangulation of $M$ to to a $\Gamma$-equivariant triangulation 
of $X$. Note that there are finitely 
isometry types of simplices in this triangulation of $X$.

The map $f'$ is defined inductively on the skeleta of the triangulation.  
On vertices we simply define $f'$ to equal $f$.  Suppose $f'$ is 
defined on the $k$-skeleton of the triangulation; then for each 
$(k+1)$-simplex $\sigma$, we have a Lipshitz map defined 
on $\partial(\sigma)$, which is a sphere.  This map extends to a lipshitz 
map on $\sigma$ by Whitney's Extension Theorem.  
Do this for each different $(k+1)$-simplex 
$\sigma$; the point is that there are only finitely many different 
lipshitz constants since there are only finitely many isometry 
types of simplices; hence the map $f'$ is lipshitz with constant the 
maximum of the Whitney lipshitz constants on the finitely many isometry 
types of simplices.
\eproof

\end{document}